\documentclass{article}

\usepackage{amssymb}

\newtheorem{thm}{Theorem}
\newtheorem{lem}{Lemma}
\newtheorem{cor}{Corollary}

\newcommand{\pf}{\noindent {\bf Proof:\,\,\,}}
\newcommand{\qed}{\hspace*{\fill}{$\rule{1ex}{1ex}$}
\medskip}

\def\a{\alpha}
\def\b{\beta}
\def\g{\gamma}

\begin{document}

\title{Hyperbolic Twistor Spaces}
\author{{D. E. Blair, J. Davidov and O. Mu\u skarov}
\footnote {The work of all three authors is supported
in part by NSF grant INT-9903302. \hfill
\hspace{1.5cm}
The second and third-named authors are members of
EDGE, Research Training Network
\hspace*{.28in}HPRN-CT-2000-00101, supported by the
European Human Potential Programme.}}

\date{}

\maketitle

\centerline{\it Dedicated to Professor Paulette
Libermann}
\vskip 20pt
\noindent
{\bf Abstract.}
In contrast to the classical twistor spaces whose
fibres are $2$-spheres, we
introduce twistor spaces over manifolds with almost
quaternionic
structures of the second kind in the sense of P.
Libermann whose fibres are
hyperbolic planes. We discuss two natural almost complex structures
on such a twistor space
and their holomorphic functions.

\vspace{0.2cm}

\vspace{0.2cm}


\vskip 20pt
\noindent
{\bf 1. Introduction}
\vskip 5pt
\noindent
In this paper we introduce hyperbolic
twistor spaces which are bundles over manifolds with almost
quaternionic structure of the second kind in the sense of
P. Libermann and whose fibres are hyperbolic planes.
These spaces admit two natural almost complex structures defined
as in the classical twistor space theory.  Our main purpose is to
study their differential-geometric properties as well as the existence
of holomorphic functions
(part of the results have been
announced in the first author's lecture
\cite{DB}).  In Section 2 we define
  hyperbolic twistor spaces and their almost
complex structures when the base manifolds are of
dimension $\geq 8$, developing the theory for
  paraquaternionic
K\"ahler manifolds.  In Section 3 we give the
corresponding treatment for base manifolds which are
4-dimensional with a metric of signature $(++--)$.
Finally in
Section 4 we treat the question of existence of
holomorphic functions on hyperbolic twistor spaces and show that,
in contrast to the classical case, there can be an abundance of (global)
holomorphic functions on a hyperbolic twistor space.

\vskip 20pt
\noindent
{\bf 2. Hyperbolic twistor spaces}
\vskip 5pt
\noindent
We begin with the
following simple observation.  In \cite{Lib} P.
Libermann introduced the
notion of an {\it almost quaternionic structure of the
second
kind} ({\it presque quaternioniennes de deuxi\`eme
esp\`ece}) on a smooth manifold $M$.
This consists of an almost complex structure $J_1$ and
an almost
product structure,
$J_2$ such that $J_1J_2+J_2J_1=0$.  Setting
$J_3=J_1J_2$ one
has a second almost product structure which also
anti-commutes with $J_1$ and $J_2$. Now on a
manifold $M$ with such a structure, set
$$j=y_1J_1+y_2J_2+y_3J_3.$$
Then $j$ is an almost complex structure on $M$ if and
only
if
$$-y_1^2+y_2^2+y_3^2=-1$$
which suggests considering a {\it hyperbolic twistor
space} $\pi:{\cal Z}\longrightarrow
M$  with fibre this hyperboloid. Recall that the
classical twistor space over a quaternionic
K\"ahler manifold is a bundle over the manifold with
the fibre being a $2$-sphere
(Salamon \cite{Sal}).

An {\it almost paraquaternionic structure} on a smooth
manifold $M$ is defined to be a rank $3$-subbundle $E$
of the endomorphisms bundle $End(TM)$ which locally is
spanned by a triple
$\{J_1,J_2,J_3\}$ which is an almost quaternionic
structure of the second kind in the sense
of P. Libermann.

There are a number of examples of almost
paraquaternionic structures including
the paraquaternionic projective space as described by
Bla\u zi\'c \cite{Bl}.  Under
certain holonomy assumptions almost paraquaternionic
structures become paraquaternionic
K\"ahler (see e.g. Garcia-Rio, Matsushita and
Vazquez-Lorenzo \cite{GMV}).
Even more strongly one has the notion of a neutral
hyperk\"ahler structure
(see Section 4) and Kamada \cite{Ka} has observed that
the only compact four-manifolds admitting such a
structure are complex tori and primary Kodaira
surfaces. We remark that the neutral hyperk\"ahler
four-manifolds are Ricci flat and self-dual
(\cite{Ka}).

The tangent bundle of
a differentiable manifold also carries an almost
paraqua\\
ternionic structure as studied by S. Ianus and C.
Udriste \cite{Ia} \cite{IU};
this includes examples where the dimension of the
manifold carrying the
structure is not necessarily $4n$.
However the most natural setting for this kind of
structure is on a manifold $M$ of dimension $4n$
with a neutral metric $g$, i.e. a pseudo-Riemannian
metric of
signature $(2n,2n)$.  One reason for this is that such
a
metric may be given with respect to which $J_1$ acts
as an
isometry on tangent spaces and $J_2$, $J_3$ act as
anti-isometries; the effect of this is that we may
define
three fundamental 2-forms $\Omega_a$, $a=1,2,3$,
by
$\Omega_a(X,Y)=g(X,J_a Y)$. If a neutral metric $g$
has this property we shall say that it is {\it
adapted} to the almost paraquaternionic structure $E$;
we shall also say that $J_1,J_2,J_3$ are {\it
compatible} with $g$.
Riemannian metrics
can be chosen such that $g(J_a X,J_a Y)=g(X,Y)$,
but then $\Omega_2$ and $\Omega_3$ are symmetric
tensor
fields instead of 2-forms.


The neutral metric $g$ induces a metric on the fibres
of
$E$ by $\displaystyle{{1\over 4n}}{\rm tr} A^tB$ where
$A$ and $B$ are
endomorphisms of $T_pM$ and $A^t$ is the adjoint of
$A$ with respect to
$g$.  This metric on the fibre is of signature
$(+--)$, the norm of
$J_1$ being $+1$ and the norms of $J_2$ and $J_3$
being
$-1$. The twistor space ${\cal Z}$ of an almost
paraquaternionic structure $E$ with an adapted neutral
metric $g$ is the unit sphere subbundle of $E$ (with
respect to the induced metric).

Alternatively one may take the Lorentz
metric $\langle\,,\rangle$ on the fibres of $E$ such
that
$\langle J_1,J_1\rangle=-1$,
$\langle J_2,J_2\rangle=+1$,
$\langle J_3,J_3\rangle=+1$.
  This metric is of signature
$(-++)$ and has the advantage of inducing
immediately a Riemannian metric of constant curvature
$-1$ on
the hyperbolic planes defined by
$-y_1^2+y_2^2+y_3^2=-1$,
in each fibre. We adopt this metric for its geometric
attractiveness but keep its negative in mind.

We will also use the following notation. For the
metric
$\langle\,,\rangle$ on the fibres of $E$ we set
$\epsilon_1=-1$ and
$\epsilon_2=\epsilon_3=+1$.  For the neutral metric
$g$ of
on the base, we set
$\varepsilon_i=\pm1$ according to the signature
$(+\cdots+-\cdots-)$.  Further, denoting
also by
$\pi$ the projection of $E$ onto $M$, if $x_i$
are local coordinates on $M$, set $q_i=x_i\circ\pi$.
We will identify the tangent space of $E$ at a point
$x\in E$ with the fibre $E_{\pi(x)}$
through that point. For a section $s$ of $E$ we denote
its vertical lift to $E$ as a vector field by $s^v$
(so $s^v=s\circ\pi$) and  frequently utilize the
natural identifications of
$J_a^v$ with $J_a$ itself and with
$\displaystyle{{\partial\over\partial y_a}}$ in terms
of the fibre
coordinates $y_1,y_2,y_3$.

An almost paraquaternionic
manifold $M$ of dimension $4n$ and neutral
metric $g$ is said to be {\it paraquaternionic
K\"ahler} if the
bundle $E$ is parallel with respect to the Levi-Civita
connection
of $g$.

As with the theory of twistor spaces over quaternionic
K\"ahler manifolds, the theory of hyperbolic twistor
spaces over
paraquaternionic K\"ahler manifolds develops nicely by
virtue of the fact that the covariant
derivatives of sections of $E$ are again sections of
$E$.  To give this development we first need
the natural machinery of horizontal lifts.

Let $D$ denote the Levi-Civita connection of the
neutral
metric on $M$.  Then the horizontal lift $X^h$ of a
vector
field $X$ to the bundle $\pi:E\longrightarrow M$ is
given by
$$X^h=\sum_iX^i{\partial\over\partial q^i}
-\sum_{a,b=1}^3
\epsilon_by_a(\langle D_XJ_a,J_b\rangle\circ\pi)
{\partial\over\partial y_b}.\eqno(2.1)$$
It is straightforward to obtain the following at a
point
$\sigma\in Z\subset E$
$$[X^h,Y^h]_\sigma=[X,Y]^h_\sigma-R(X,Y)\sigma$$
(we adopt here the following definition of the
curvature tensor:
$R(X,Y)=\nabla_{X}\nabla_{Y}-\nabla_{Y}\nabla_{X}-\nabla_{[X,Y]}$).
For a section $s$ of $E$,
$$[X^h,s^v]=(D_Xs)^v.$$

\noindent
Define a metric on $E$ by $h_t=\pi^*g
+t\langle\,,\rangle$, $t\neq 0$,
and for simplicity denote its Levi-Civita connection
by $\bar\nabla$ instead of
$\bar\nabla^t$, but note the $t$ in the formulas
below.  Then
$$(\bar\nabla_{X^h}Y^h)_\sigma=(D_XY)^h_\sigma
-{1\over 2}R(X,Y)\sigma,\eqno(2.2)$$
$$(\bar\nabla_{X^h}s^v)_\sigma
={t\over 2}(\hat R_{\sigma
s}X)^h+(D_Xs)^v_\sigma\eqno(2.3)$$
where $g(\hat R_{\sigma s}X,Y)=\langle
R(X,Y)\sigma,s^v\rangle$ and
$$(\bar\nabla_{s^v}X^h)_\sigma
={t\over 2}(\hat R_{\sigma s}X)^h,\quad
\bar\nabla_{k^{v}}s^{v}=0.\eqno(2.4)$$
for sections $k$ and $s$ of $E$.

\begin{lem}
The Weingarten map $A_t$ of the hyperbolic twistor
space $Z$, as a
hypersurface in the bundle space $E$,
annihilates horizontal vectors and acts on vertical
vectors by
$A_tV=\displaystyle{\sqrt{|t|}\over t}V$.
\end{lem}
\pf
The position vector $\sigma$ gives rise to a normal
$\displaystyle\nu_t={\sigma\over\sqrt{|t|}}$ to the
hyperboloid in
the fibres of $E$
and the Weingarten map is given by
$A_tX=\displaystyle{|t|\over t}\bar\nabla_X\nu_t$.
Using the summation convention for repeated indices
$a,b,c=1,2,3$, we have
first for a horizontal lift (equation (2.1)) to a
point $\sigma$
$$\bar\nabla_{X^h}\nu_t={1\over\sqrt{|t|}}\bar\nabla_{X^h}y_c{\partial\over\partial
y_c}$$
$$={1\over\sqrt{|t|}}\Big(-\epsilon_cy_a(\langle
D_XJ_a,J_c\rangle\circ\pi)
{\partial\over\partial y_c}\Big)
+{1\over\sqrt{|t|}}y_c\big({t\over 2}(\hat R_{\sigma
{J_c}}X)^h+(D_XJ_c)_\sigma\big)$$
$$={t\over 2\sqrt{|t|}}(\hat R_{\sigma\sigma}X)^h=0.$$
Similarly for a vertical tangent vector $V$
$$\bar\nabla_V\nu_t={1\over\sqrt{|t|}}\bar\nabla_Vy_c{\partial\over\partial
y_c}
={1\over\sqrt{|t|}}(Vy_c){\partial\over\partial y_c}
={V^c\over\sqrt{|t|}}{\partial\over\partial
y_c}={1\over\sqrt{|t|}}V.$$
Thus $A_tV=\displaystyle{\sqrt{|t|}\over t}V$ for a
vertical tangent
vector $V$ and $A_tX=0$ for a
horizontal tangent vector $X$.
\qed

    We now define two almost complex structures ${\cal
J}_1$ and
${\cal J}_2$ on the hyperbolic twistor space ${\cal
Z}$ as follows. Acting on
horizontal vectors these are the same and given by
${\cal J}_1X^h_\sigma={\cal
J}_2X^h_\sigma=(jX)^h_\sigma$ where as before
$j=\sum y_aJ_a$ is the point ${\sigma}$ considered as
an endomorphism of $TM$.
For a vertical vector
$V=\displaystyle{V^1{\partial\over\partial y_1}
+V^2{\partial\over\partial y_2}
+V^3{\partial\over\partial y_3}}$ tangent to
${\cal Z}$, i.e.
$\langle\sigma,V\rangle=0$, let
$${\cal J}_1V=
(y_3V^2-y_2V^3){\partial\over\partial y_1}
+(y_3V^1-y_1V^3){\partial\over\partial y_2}
+(y_1V^2-y_2V^1){\partial\over\partial
y_3}\eqno(2.5)$$
and let ${\cal J}_2V$ be the negative of this
expression.
At each point $p\in M$, the local endomorphisms
$\{J_1,J_2,J_3\}$ form an orthonormal basis of the
fibre $E_p$ of $E$ over $p$ and define the same
orientation on it. Denote by $\times$ the vector
product on the $3$-dimensional vector space $E_p$
determined by this orienation and the metric of $E$
(in other words determined by the paraquaternionic
algebra). Then
${\cal J}_kV=(-1)^{k-1}\sigma\times V$, $k=1,2$, for
any $\sigma\in {\cal Z}$

  Define a pseudo-Riemannian metric on ${\cal Z}$ by
$h_t=\pi^*g+t\langle\,,\rangle_v$, $t\neq 0$,
$\langle\,,\rangle_v$ being the restriction of
$\langle\,,\rangle$
to the fibres (hyperbolic planes) of ${\cal Z}$ and
denote the Levi-Civita connection of $h_t$
by $\nabla$ for simplicity. It is easy to check that
this metric is Hermitian with respect to both ${\cal
J}_1$ and
${\cal J}_2$.

We now review paraquaternionic
K\"ahler geometry;
the development of theory of
paraquaternionic K\"ahler structures was
carried out by
Garcia-Rio, Matsushita and Vazquez-Lorenzo \cite{GMV}
when the dimension
of the base manifold is $4n\geq 8$.  As with the
theory of quaternionic K\"ahler manifolds, dimension 4
is
special.  At the beginning however we will retain the
4-dimensional case and point out where the differences
occur.
The parallel to the present development in the
quaternionic
K\"ahler case can be found in Ishihara \cite{Is}.

  Let $\{J_1,J_2,J_3\}$ be a local almost quaternionic
structure of the second kind which
spans the bundle $E$. Since $E$ is parallel with
respect to $D$, there exist local 1-forms
$\a$, $\b$ and $\g$ such that

$$D_XJ_1=\hskip 1cm -\g(X)J_2-\b(X)J_3,$$
$$D_XJ_2=-\g(X)J_1\hskip 1cm-\a(X)J_3,\eqno(2.6)$$
$$D_XJ_3=-\b(X)J_1+\a(X)J_2.\hskip 1cm$$

 From the group theoretic point of view, this structure
corresponds
to the linear holonomy group being a subgroup of
$Sp(n,{\Bbb R})\cdot Sp(1,{\Bbb R})$, just as a
quaternionic
K\"ahler  structure corresponds
to the linear holonomy group being a subgroup of
$Sp(n)\cdot Sp(1)$.  For $n=1$ this is not a
restriction.

Setting
$$A=2(d\a-\b\wedge\g),\quad B=2(d\b-\a\wedge\g),
\quad C=2(d\g+\a\wedge\b)$$
one can easily obtain the following central relation
of the
action of the curvature tensor:

$$R(X,Y)J_1=\hskip 1cm -C(X,Y)J_2-B(X,Y)J_3,$$
$$R(X,Y)J_2=-C(X,Y)J_1\hskip
1cm-A(X,Y)J_3,\eqno(2.7)$$
$$R(X,Y)J_3=-B(X,Y)J_1+A(X,Y)J_2.\hskip 1cm$$

Moreover a paraquaternionic K\"ahler manifold of
dimension $\geq 8$ is Einstein and
$A$, $B$, $C$ satisfy
$$A(X,Y)=-{\tau g(X,J_1Y)\over 4n(n+2)},\;\,
B(X,Y)=-{\tau g(X,J_2Y)\over 4n(n+2)},\;\,
C(X,Y)={\tau g(X,J_3Y)\over 4n(n+2)}\eqno(2.8)$$
where $\tau$ is the scalar curvature of the metric
$g$.

We now give our first result.

\begin{thm}
On the hyperbolic twistor space of a paraquaternionic
K\"ahler manifold of
dimension $4n\geq 8$ we have the following:
\begin{enumerate}
\item[$(i)$] The almost complex structure ${\cal J}_1$
is integrable and the Hermitian structure
$({\cal J}_1,h_t)$ is semi-K\"ahler. It is indefinite
K\"ahler if and only if $t\tau=-4n(n+2)$.
\item[$(ii)$]The almost complex structure ${\cal J}_2$
is never integrable but the almost Hermitian structure
$({\cal J}_2,h_t)$ is semi-K\"ahler. It is indefinite
almost K\"ahler  if and only if $t\tau=4n(n+2)$ and
indefinite nearly K\"ahler if and only if
$t\tau=-2n(n+2)$.
\end{enumerate}
\end{thm}

\pf
The major effort of the proof is to compute the
covariant derivatives of ${\cal J}_i$,
$i=1,2$. To begin, by Lemma 1 and (2.2)
$$(\nabla_{X^h}{\cal J}_i)Y^h\big|_\sigma
=\bar\nabla_{X^h}{\cal J}_iY^h-{\cal
J}_i\bar\nabla_{X^h}Y^h$$
$$=\bar\nabla_{X^h}(y_1J_1Y+y_2J_2Y+y_3J_3Y)^h-{\cal
J}_i(D_XY)^h
+{1\over 2}{\cal J}_iR(X,Y)\sigma$$
Expanding furthur and using (2.1) and (2.7),
$$(\nabla_{X^h}{\cal J}_i)Y^h\big|_\sigma
=-{1\over 2}y_1R(X,J_1Y)\sigma
-{1\over 2}y_2R(X,J_2Y)\sigma
-{1\over 2}y_3R(X,J_3Y)\sigma$$
$$+{1\over 2}{\cal
J}_i\Big[y_1\Big(-C(X,Y){\partial\over\partial
y_2}-B(X,Y){\partial\over\partial y_3}\Big)
+y_2\Big(-C(X,Y){\partial\over\partial
y_1}-A(X,Y){\partial\over\partial y_3}\Big)$$
$$+y_3\Big(-B(X,Y){\partial\over\partial
y_1}+A(X,Y){\partial\over\partial y_2}\Big)\Big].$$
Note that $\displaystyle{y_1{\partial\over\partial
y_2}+y_2{\partial\over\partial y_1}}$, etc. are
tangent
to the fibres.  Applying ${\cal J}_i$ and expanding
the curvature terms by (2.7) we have
$$(\nabla_{X^h}{\cal
J}_1)Y^h\big|_\sigma=0\eqno(2.9)$$
and
$$(\nabla_{X^h}{\cal J}_2)Y^h\big|_\sigma
=-\Big[\big((y_2^2+y_3^2)A(X,Y)+y_1y_2B(X,Y)-y_1y_3C(X,Y)\big){\partial\over\partial
y_1}$$
$$+\big((y_1y_2A(X,Y)+(y_1^2-y_3^2)B(X,Y)-y_2y_3C(X,Y)\big){\partial\over\partial
y_2}$$
$$+\big((y_1y_3A(X,Y)+y_2y_3B(X,Y)+(y_2^2-y_1^2)C(X,Y)\big){\partial\over\partial
y_3}\Big]$$
$$={\tau\over 4n(n+2)}
\Big[\big(-g(X,J_1Y)+y_1g(X,jY)\big){\partial\over\partial
y_1}$$
$$+\big(g(X,J_2Y)+y_2g(X,jY)\big){\partial\over\partial
y_2}
+\big(g(X,J_3Y)+y_3g(X,jY)\big){\partial\over\partial
y_3}\Big]$$
using (2.8) and $-y_1^2+y_2^2+y_3^2=-1$.  Taking the
inner product with a vertical tangent
vector $V$ we have
$$h_t((\nabla_{X^h}{\cal J}_2)Y^h,V)
={{t\tau}\over
4n(n+2)}\big[V^1g(X,J_1Y)+V^2g(X,J_2Y)+V^3g(X,J_3Y)\big].\eqno(2.10)$$

    For $(\nabla_{X^h}{\cal J}_i)V$, its horizontal
part may be found immediately from the above
and to show that its vertical part vanishes we show
that $(\nabla_{X^h}{\cal J}_i)V$ is
horizontal.  To do this effectively recall that ${\cal
J}_1$ can be given by equation (2.5)
and we may regard this formula as extended to $E$,
i.e. ${\cal J}_1V$ is given by this
formula for $V$ tangent to $E$, even though one no
longer has ${\cal J}_1^2=-I$.  Then
$$(\nabla_{X^h}{\cal J}_1){\partial\over\partial y_1}
=\bar\nabla_{X^h}\Big(y_3{\partial\over\partial
y_2}-y_2{\partial\over\partial y_3}\Big)
-{\cal J}_1\Big({t\over 2}(\hat R_{\sigma {J_1}}X)^h
-\g(X){\partial\over\partial
y_2}-\b(X){\partial\over\partial y_3}\Big)$$
$$=\big(y_1\b(X)+y_2\a(X)\big){\partial\over\partial
y_2}
+y_3\Big({t\over 2}(\hat R_{\sigma {J_2}}X)^h
-\g(X){\partial\over\partial
y_1}-\a(X){\partial\over\partial y_3}\Big)$$
$$+\big(-y_1\g(X)+y_3\a(X)\big){\partial\over\partial
y_3}
-y_2\Big({t\over 2}(\hat R_{\sigma {J_3}}X)^h
-\b(X){\partial\over\partial
y_1}+\a(X){\partial\over\partial y_2}\Big)$$
$$-{t\over 2}(j\hat R_{\sigma {J_1}}X)^h
+\g(X)\Big(y_3{\partial\over\partial
y_1}+y_1{\partial\over\partial y_3}\Big)
+\b(X)\Big(-y_2{\partial\over\partial
y_1}-y_1{\partial\over\partial y_2}\Big)$$
$$={t\over 2}y_3(\hat R_{\sigma {J_2}}X)^h
-{t\over 2}y_2(\hat R_{\sigma {J_3}}X)^h-{t\over
2}(j\hat R_{\sigma {J_1}}X)^h$$
which is horizontal.  The proof for
$\displaystyle{{\partial\over\partial y_2}}$ and
$\displaystyle{{\partial\over\partial y_3}}$ and
for ${\cal J}_2$ is similar.

Similarly treating $(\nabla_V{\cal J}_i)X^h$ we find
that
$$h_t((\nabla_V{\cal J}_i)X^h,Y^h)$$
$$={4n(n+2)+t\tau\over
4n(n+2)}\big[V^1g(J_1X,Y)+V^2g(J_2X,Y)+V^3g(J_3X,Y)\big].
\eqno(2.11)$$

Finally for vertical tangent vectors $V$ and $W$
$$(\nabla_V{\cal J}_i)W=\bar\nabla_V{\cal
J}_iW-t\langle A_tV,{\cal J}_iW\rangle\nu_t
-{\cal J}_i\bar\nabla_VW$$
noting that the extension of ${\cal J}_i$ to $E$
annihilates $\nu_t$.  Treating the terms
separately for ${\cal J}_1$
$$\bar\nabla_V{\cal
J}_1W=(V^3W^2+y_3VW^2-V^2W^3-y_2VW^3){\partial\over\partial
y_1}$$
$$+(V^3W^1+y_3VW^1-V^1W^3-y_1VW^3){\partial\over\partial
y_2}$$
$$+(V^1W^2+y_1VW^2-V^2W^1-y_2VW^1){\partial\over\partial
y_3},$$
$${\cal J}_1\bar\nabla_VW
=(VW^1)\Big(y_3{\partial\over\partial
y_2}-y_2{\partial\over\partial y_3}\Big)
+(VW^2)\Big(y_3{\partial\over\partial
y_1}+y_1{\partial\over\partial y_3}\Big)$$
$$+(VW^3)\Big(-y_2{\partial\over\partial
y_1}-y_1{\partial\over\partial y_2}\Big)$$
and using
$\langle\sigma,V\rangle=\langle\sigma,W\rangle=0$ and
$-y_1^2+y_2^2+y_3^2=-1$
$$y_1\langle V,{\cal J}_1W\rangle=V^3W^2-V^2W^3,$$
$$y_2\langle V,{\cal J}_1W\rangle=V^3W^1-V^1W^3,$$
$$y_3\langle V,{\cal J}_1W\rangle=V^1W^2-V^2W^1.$$
Combining these we have $(\nabla_V{\cal J}_1)W=0$ and
similarly $(\nabla_V{\cal J}_2)W=0$.

Using these computations we can now easily complete
the proof of
Theorem 1.  That the almost Hermitian structure
  $({\cal J}_1,h)$ is K\"ahler if and only if
$t\tau=-4n(n+2)$ follows immediately from the above
relations, especially equations (2.9) and
(2.11).  To show the integrability of
  ${\cal J}_1$ first recall that the Nijenhuis tensor
$N_i$ of ${\cal J}_i$
$$N_i(X,Y)=-[X,Y]+[{\cal J}_iX,{\cal J}_iY]-{\cal
J}_i[{\cal J}_iX,Y]-{\cal J}_i[X,{\cal J}_iY]$$
  may be written in terms
of the connection $\nabla$ as
$$N_i(X,Y)
={\cal J}_i(\nabla_Y{\cal J}_i)X-(\nabla_{{\cal
J}_iY}{\cal J}_i)X
-{\cal J}_i(\nabla_X{\cal J}_i)Y+(\nabla_{{\cal
J}_iX}{\cal J}_i)Y.\eqno (2.12)$$
The cases $N_1(X^h,Y^h)=0$ and $N_1(V,W)=0$ are
immediate. For $N_1(V,X^h)$, observe that the first
two terms of the
expansion (2.12) vanish while the remaining two are
horizontal.  Thus it is enough to compute
$$h_t(N_1(V,X^h),Y^h)
=h_t((\nabla_V{\cal
J}_1)X^h,(jY)^h)+h_t((\nabla_{{\cal J}_1V}{\cal
J}_1)X^H,Y^h);$$
upon expansion using (2.11) the two terms will cancel.

For the almost Hermitian structure
  $({\cal J}_2,h)$, to see that it is almost K\"ahler
if and only if
$t\tau=4n(n+2)$, the key case to consider is
$$h_t((\nabla_{X^h}{\cal
J}_2)Y^h,V)+h_t((\nabla_V{\cal J}_2)X^h,Y^h)
+h_t((\nabla_{Y^h}{\cal J}_2)V,X^h)$$
$$={t\tau-4n(n+2)\over
4n(n+2)}\big[V^1g(X,J_1Y)+V^2g(X,J_2Y)+V^3g(X,J_3Y)\big].$$
To see that $({\cal J}_2,h)$ is nearly K\"ahler if and
only if
$t\tau=-2n(n+2)$, note that
$h_t((\nabla_{X^h}{\cal J}_2)Y^h+(\nabla_{Y^h}{\cal
J}_2)X^h,V)=0$ by the skew-symmetry
in equation (2.10) and by equations (2.10) and (2.11)
$$h_t((\nabla_{X^h}{\cal J}_2)V+(\nabla_V{\cal
J}_2)X^h,Y^h)$$
$$={4n(n+2)+2t\tau\over
4n(n+2)}\big[V^1g(J_1X,Y)+V^2g(J_2X,Y)+V^3g(J_3X,Y)\big].$$

   To show the non-integrability of ${\cal J}_2$, we
compute
$h_t(N_2(V,X^h),Y^h)$ at the point $(1,0,0)$ with
$V=\displaystyle{{\partial\over\partial
y_2}+{\partial\over\partial y_3}}$.
The first term in the expansion (2.12) yields
$$-h_t((\nabla_{X^h}{\cal J}_2)V,(jY)^h)
={t\tau\over 4n(n+2)}\big[-g(X,J_3Y)+g(X,J_2Y)\big].$$
Proceeding in this way with the other terms we get
$$h_t(N_2(V,X^h),Y^h)=-2[g(X,J_2Y)-g(X,J_3Y)]$$
which is not identically zero, e.g take $X=J_2Y$.

   Finally, note that, by the above computations,
$\displaystyle{(\nabla_{X^h}{\cal J}_{i})X^h=0}$ for
any $X\in TM$ and $\displaystyle{(\nabla_V{\cal
J}_{i})V=0}$ for any vertical vector $V$, $i=1,2$.
This implies that $({\cal J}_i,h_t)$ has co-closed
fundamental $2$-form, i.e. it is semi-K\"ahler.
\qed

\noindent
{\bf Remark.}
The values of the scalar curvature appearing in
Theorem 1 for $t=1$ are the negatives of what
one has in the usual twistor space over a quaternionic
K\"ahler manifold of dimension $\geq 8$, see e.g.
\cite{AGI}.
This sign change is due to
our choice of metric on the fibres of $E$.  If we take
$<\,,>$ as the
$(+--)$ metric we would have the other values, but the
fibres of ${\cal Z}$
would then have a negative definite metric.
  In the classical case the almost complex structure
${\cal J}_1$ was
introduced and shown to be integrable
by S. Salamon \cite{Sal} and independently by L.
B\'erard Bergery (unpublished but see e.g.
Besse
\cite{bes}).

\vskip 20pt

\noindent
{\bf 3. The 4-dimensional case}
\vskip 5pt
\noindent
   Let $M$ be an oriented 4-dimensional manifold with a
neutral
metric $g$ and ${\bf e}_1,\ldots,{\bf e}_4$ a local
orthonormal frame
with
${\bf e}_1\wedge{\bf e}_2\wedge{\bf e}_3\wedge{\bf
e}_4$ giving
the orientation.
The metric $g$ induces a metric on bundle of
bivectors,
$\bigwedge^2 TM$, by
$$g({\bf e}_i\wedge{\bf e}_j,{\bf e}_k\wedge{\bf e}_l)
={1\over 2}\left|\begin{array}{cc}
\varepsilon_i\delta_{ik}&\varepsilon_i\delta_{il}\\
\varepsilon_j\delta_{jk}&\varepsilon_j\delta_{jl}
\end{array}\right|,\quad
\varepsilon_1=\varepsilon_2=1,\quad
\varepsilon_3=\varepsilon_4=-1.$$
The Hodge star operator of the neutral metric acting
on
$\bigwedge^2 TM$ is given by
$$*({\bf e}_1\wedge{\bf e}_2)={\bf e}_3\wedge{\bf
e}_4,\quad
*({\bf e}_1\wedge{\bf e}_3)={\bf e}_2\wedge{\bf
e}_4,\quad
*({\bf e}_1\wedge{\bf e}_4)=-{\bf e}_2\wedge{\bf
e}_3.$$
Let $\bigwedge^-$ and $\bigwedge^+$ denote the
subbundles of
$\bigwedge^2 TM$ determined by the corresponding
eigenvalues of
the Hodge star operator. The metrics induced on
$\bigwedge^-$  and
$\bigwedge^+$ have signature $(+--)$.

  Setting
$$s_1={\bf e}_1\wedge{\bf e}_2-{\bf e}_3\wedge{\bf
e}_4,\quad\quad
\bar s_1={\bf e}_1\wedge{\bf e}_2+{\bf e}_3\wedge{\bf
e}_4,$$
$$s_2={\bf e}_1\wedge{\bf e}_3-{\bf e}_2\wedge{\bf
e}_4,\quad\quad
\bar s_2={\bf e}_1\wedge{\bf e}_3+{\bf e}_2\wedge{\bf
e}_4,$$
$$s_3={\bf e}_1\wedge{\bf e}_4+{\bf e}_2\wedge{\bf
e}_3,\quad\quad
\bar s_3={\bf e}_1\wedge{\bf e}_4-{\bf e}_2\wedge{\bf
e}_3,$$
$\{s_1,s_2,s_3\}$ and $\{\bar s_1,\bar s_2,\bar s_3\}$
are local
oriented orthonormal frames for $\bigwedge^-$ and
$\bigwedge^+$
respectively.

   Reversing the orientation of $M$ just interchanges
the roles of
$\bigwedge^-$ and $\bigwedge^+$, and we shall
concentrate only on the
bundle $\bigwedge^-$.

   Further we shall often identify $\bigwedge^2 TM$
with the bundle of skew-symmetric endomorphisms of
$TM$ by the correspondance that assigns
to each $\sigma\in\bigwedge^2 TM$ the endomorphism
$K_{\sigma}$ on $T_pM$, $p=\pi(\sigma)$, defined
by
$$g(K_{\sigma}X,Y)=2g(\sigma,X\wedge Y); X,Y\in
T_pM.\eqno(3.1)$$
Thus $s_1,s_2,s_3$ are identified with the
endomorphisms representing in the frame ${\bf
e}_1,\ldots,{\bf e}_4$ by the matricies
$$\left( \begin{array}{cccc}
           0 &-1 & 0 & 0\\
           1 & 0 & 0 & 0\\
           0 & 0 & 0 &-1 \\
                0 & 0 & 1 & 0
        \end{array} \right),\;
\left( \begin{array}{cccc}
          0 &  0 &1 & 0\\
          0 & 0 & 0 &-1\\
         1 & 0 & 0 & 0 \\
                0 &-1 & 0 & 0
        \end{array} \right),\;
\left( \begin{array}{cccc}
           0 & 0 & 0 & 1\\
          0 & 0 & 1 & 0\\
           0 & 1 & 0 & 0 \\
                1 & 0 & 0 & 0
        \end{array} \right)$$
Hence the bundle $E=\bigwedge^-$ defines an almost
paraquaternionic structure on $M$, the local
endomorphisms $\{J_1,J_2,J_3\}$ spanning $E$ being
$J_1=K_{s_1}, J_2=K_{s_2}, J_3=K_{s_3}$. Moreover, the
Levi-Civita connection of $M$ preserves the bundle
$\bigwedge^-$. So, as we have already mentioned,
the existence of a paraquaternionic K\"ahler structure
does not impose any restriction on the oriented
Riemannian four-manifolds (the four-dimensional analog
of paraquaternionic K\"ahler manifolds are the
Einstein self-dual manifolds).

    Now, in accordance with Section 2, the hyperbolic
twistor space ${\cal Z}$ of $M$ is defined to be the
unit sphere bundle in $\bigwedge^-$. It can be
identified via (3.1) with the space of all complex
structures on the tangent spaces of $M$ compatible
with its metric and orientation. We keep the notations
${\cal J}_1$
and ${\cal J}_2$ for the natural almost complex
structures on ${\cal Z}$ noting that in the Riemannian
case they have been introduced and studied by
Atiyah-Hitchin-Singer \cite{AHS} and, respectively,
Eells-Salamon  \cite{ES}.

   Let ${\cal R}:\bigwedge^2
TM\longrightarrow\bigwedge^2 TM$ be the curvature
operator of $(M,g)$. It is related to the curvature
tensor $R$ by
$$g({\cal R}(X\wedge Y),Z\wedge T)=-g(R(X,Y)Z,T);
X,Y,Z,T\in TM.$$
It is not hard to check that, for any $a\in\bigwedge^2
TM$ and $b,c\in\bigwedge^-$, we have
$$g(R(a)b,c)=-g({\cal R}(b\times c),a)\eqno
(3.2)$$
where $R$ on the left-hand side stands for the
curvature of the connection on the
bundle $\bigwedge^2 TM$ induced by the Levi-Civita
connection of $M$.

    Let us also note that if $V\in {\cal V}_{\sigma}$
and $X,Y\in T_{\pi(\sigma)}M$, then
$$g(\sigma\times V,X\wedge Y)=-g(V,X\wedge
K_{\sigma}Y).\eqno
(3.3)$$

    The curvature operator ${\cal
R}:\bigwedge^2TM\longrightarrow\bigwedge^2 TM$
admits an $SO(2,2)$-irreducible decomposition
$${\cal R}={\tau\over 6}I+{\cal B}+{\cal W}^++{\cal
W}^-$$
similar to that in the 4-dimensional Riemannian case.
Here $\cal B$ represents the
traceless Ricci tensor,  ${\cal W}={\cal W}^++{\cal
W}^-$ corresponds to the Weyl conformal tensor, and
${\cal W}^{\pm}=\displaystyle{\frac{1}{2}({\cal
W}\pm\ast{\cal W})}$.  The metric
$g$ is said to be {\it self-dual} if ${\cal W}^-=0$.

    The metric $g$ on the bundle $\pi:\bigwedge^2
TM\longrightarrow M$ induced by the metric of $M$ is
negative definite on the fibres of ${\cal Z}$ and, as
in Section 2, we adopt the metric
$\langle\,,\rangle=-g$ on $\bigwedge^2 TM$. Setting
$h_t=\pi^*g +t\langle\,,\rangle$
for any real $t\neq 0$ we get a $1$-parameter family
of pseudo-Riemannian metrics on ${\cal Z}$
compatible with the almost complex structures ${\cal
J}_1$ and ${\cal J}_2$

     Again for simplicity we denote by $\nabla$ the
Levi-Civita connection of $({\cal Z},h_t)$ and let
$D$ be the Levi-Civita connection of $(M,g)$.

    Let $X,Y$ be vector fields on $M$ and $V$ a
vertical vector field on ${\cal Z}$. Then, for any
point $\sigma\in {\cal Z}$,
$$(\nabla_{X^h}Y^h)_{\sigma}=(D_{X}Y)_{\sigma}^h-{1\over
2}R(X,Y)\sigma,\eqno
(3.4)$$
$$(\nabla_{V}X^h)_{\sigma}={\cal
H}(\nabla_{X^h}V)_{\sigma}=
-{t\over 2}(R(\sigma\times V)X)_{\sigma}^h.\eqno
(3.5)$$
where ${\cal H}$ means "the horizontal component".

    Indeed, the first identity is a consequence of the
standard formula for the Levi-Civita connection and
the fact that
$[X^h,Y^h]_{\sigma}=[X,Y]_{\sigma}^h-R(X,Y)\sigma$,
$\sigma\in {\cal Z}$. To see (3.5), let us note that
$\nabla_{V}X^h$ is a horizontal vector field since the
fibres of ${\cal Z}$ are totally geodesic
submanifolds. On the other hand, $[V,X^h]$ is a
vertical vector field, hence $\nabla_{V}X^h=
{\cal H}\nabla_{X^h}V$. Then, by (3.2), we have
$$h_t(\nabla_{V}X^h,Y^h)=h_t(\nabla_{X^h}V,Y^h)=-h_t(V,\nabla_{X^h}Y^h)=
-{t\over 2}g(R(X,Y)\sigma,V)=
$$
$$
{t\over 2}g({\cal R}(\sigma\times V),X\wedge Y)=
-{t\over 2}h_t((R(\sigma\times
V)X)^h_{\sigma},Y^h_{\sigma})$$
and we get the second identity in (3.5).

     We are now going to compute the covariant
derivative $\nabla{\cal J}_k$ of the almost complex
structure ${\cal J}_k$ on the twistor space ${\cal
Z}$, $k=1,2$. The computation is similar to that in
\cite{G, M} and we present it here for completeness.

   Let $\Omega_{k,t}(A,B)=h_t(A,{\cal J}_{k}B)$ be the
fundamental 2-form of the almost Hermitian structure
$({\cal J}_k,h_t)$.

\begin{lem}
  Let $\sigma\in {\cal Z}, X,Y\in T_{\pi(\sigma)}M$ and
$V\in {\cal V}_{\sigma}$. Then:
$$
(\nabla_{X^h}\Omega_{k,t})(Y^h,V)_{\sigma}={t\over
2}[(-1)^{k}g({\cal R}(V),X\wedge Y)+ g({\cal
R}(\sigma\times V),X\wedge K_{\sigma}Y)]
$$
$$
(\nabla_{V}\Omega_{k,t})(X^h,Y^h)_{\sigma}=-2g(V,X\wedge
Y) - {t\over 2}g({\cal R}(\sigma\times V),X\wedge
K_{\sigma}Y+K_{\sigma}X\wedge Y)
$$
$$
(\nabla_{A}\Omega_{k,t})(B,C)=0
$$
when $A,B,C$ are horizontal vectors or at least two of
them are vertical.
\end{lem}
\pf
Extend $X,Y$ to vector fields in a neighborhood of the
point $p=\pi(\sigma)$.
Then, by (3.4), (3.5) and (3.2), we have
$$
(\nabla_{X^h}\Omega_{k,t})(Y^h,V)=-h_t(\nabla_{X^h}Y^h,{\cal
J}_kV)+h_t({\cal J}_kY^h,\nabla_{X^h}V)
$$

$$
=(-1)^{k}{t\over 2}g(R(X,Y)\sigma,\sigma\times
V)+h_t((K_{\sigma}Y)^h,[X^h,V]+\nabla_{V}X^h)
$$

$$
=(-1)^{k}{t\over 2}g({\cal R}(V),X\wedge Y) + {t\over
2}g({\cal R}(\sigma\times V),X\wedge K_{\sigma}Y).
$$

    Next, by (3.5), we have
$$(\nabla_{V}\Omega_{k,t})(X^h,Y^h)=Vh_t(X^h,{\cal
J}_kY^h)-{t\over 2}g({\cal R}(\sigma\times V),X\wedge
K_{\sigma}Y+K_{\sigma}X\wedge Y)$$
Moreover, $h_t(X^h,{\cal
J}_kY^h)=2\sum_{a=1}^3y_a(g(s_a,Y\wedge X)\circ\pi),$
hence $Vh_t(X^h,{\cal J}_kY^h)=g(V,Y\wedge X)$.

   Let $U,V,W$ be vertical vector fields on ${\cal Z}$
near the point $\sigma$. Then
$$(\nabla_{U}\Omega_{k,t})(V,W)=0$$
since the fibres of ${\cal Z}$ are totally geodesic
submanifolds and the restriction of ${\cal J}_k$ on
each fibre is K\"ahlerian. We also have
$(\nabla_{U}\Omega_{k,t})(X^h,V)=0$ in view of (3.5)
and the fact that $\nabla_{U}V$ is a vertical vector
field. Next, by (3.5), we have:
$$
(\nabla_{X^h}\Omega_{k,t})(V,W)=h_t(V,\nabla_{X^h}{\cal
J}_kW-{\cal J}_k\nabla_{X^h}W)
$$
$$
=h_t(V,[X^h,{\cal J}_kW]-{\cal J}_k[X^h,W])=0
$$
since, as is easy to see, $[X^h,{\cal J}_kW]={\cal
J}_k[X^h,W]$. Indeed, take an oriented orthonormal
frame
${\bf e}_1,\ldots,{\bf e}_4$ of $TM$ near the point
$p=\pi(\sigma)$ such that $D{\bf e}_i|_p=0$,
$1\leq i\leq 4$, and $s_1(p)=\sigma$ ($s_1,s_2,s_3$
are defined by means of ${\bf e}_1,\ldots{\bf e}_4$ as
in the beginning of this section). The vector fields
$$
U=y_2{\partial\over\partial
y_1}+y_1{\partial\over\partial y_2},\quad
J_1U=y_1y_3{\partial\over\partial
y_1}+y_2y_3{\partial\over\partial y_2}+
(1+y_3^2){\partial\over\partial y_3}
$$
form a frame for the vertical bundle on ${\cal Z}$
near $\sigma$. Since $Ds_k|_p=0$, we have
$[X^h,U]_{\sigma}=[X^h,{\cal J}_1U]_{\sigma}=0$. It
follows that for every vertical vector field $W$,
$[X^h,{\cal J}_kW]_{\sigma}={\cal
J}_k[X^h,W]_{\sigma}$.

     Finally, let $X,Y,Z$ be vector fields on $M$ near
$p$. Then
$$(\nabla_{X^h}\Omega_{k,t})(Y^h,Z^h)_{\sigma}=-g(D_{X}s_1,Y\wedge
Z)_p=0$$
since $D_{X}s_1|_p=0$.
\qed

    Recall that
$$h_t(N_k(A,B),C)=-(\nabla_{A}\Omega_{k,t})(B,{\cal
J}_kC)+
(\nabla_{B}\Omega_{k,t})(A,{\cal J}_kC)
$$
$$
-(\nabla_{{\cal J}_kA}\Omega_{k,t})(B,C)
+(\nabla_{{\cal J}_kB}\Omega_{k,t})(A,C)$$
  where $N_k$ is the Nijenhuis tensor of the almost
complex structure ${\cal J}_k$. Then Lemma 2, (3.2)
and (3.3) imply the following:

\begin{cor}
Let $\sigma\in {\cal Z}$, $X,Y\in T_{\pi(\sigma})M$
and $V,W\in {\cal V}_{\sigma}$. Then:
$$
N_k(X^h,Y^h)_{\sigma}=R(X\wedge Y-K_{\sigma}X\wedge
K_{\sigma}Y)\sigma +
(-1)^{k-1}\sigma\times R(K_{\sigma}X\wedge Y+X\wedge
K_{\sigma}Y)\sigma,
$$
$$
N_k(X^h,V)_{\sigma}=2[-1+(-1)^{k-1}]g(V,X\wedge
K_{\sigma}Y),
$$

$$N_k(V,W)=0.$$
\end{cor}

\begin{cor}
  Let $A,B,C\in T_{\sigma}{\cal Z}$ and set
$X=\pi_{\ast}A, Y=\pi_{\ast}B, Z=\pi_{\ast}C$ and
$U={\cal V}_{\sigma}A, V={\cal V}_{\sigma}B, W={\cal
V}_{\sigma}C$. Then:
$$
3d\Omega_{k,t}(A,B,C)=t(-1)^{k}[g(R(U),Y\wedge
Z)+g(R(V),Z\wedge X)+
                                     g(R(W),X\wedge Y)]
$$
$$
                         - 2[g(U,Y\wedge Z)+g(V,Z\wedge
X)+g(W,X\wedge Z)].
$$


\end{cor}

\begin{cor}
   Let $A\in T_{\sigma}{\cal Z}$ and $U={\cal
V}_{\sigma}A$. Then the co-differential of
$\Omega_{k,t}$
is given by
$$
\delta\Omega_{k,t}(A)=tg(R(\sigma)\sigma,U).
$$
\end{cor}

\begin{thm}
On the hyperbolic twistor space ${\cal Z}$ of an
oriented $4$-dimensional manifold $M$
with a neutral metric $g$ we have the following:
\begin{enumerate}
\item[$(i)$] The almost complex structure ${\cal J}_1$
is integrable if and only if the metric $g$ is
self-dual. The almost Hermitian structure $({\cal
J}_1,h_t)$ is semi-K\"ahler if and only if $g$ is
self-dual.
It is indefinite K\"ahler if and
only if the metric $g$ is Einstein, self-dual, and
$t\tau=-12$.
\item[$(ii$)] The almost complex structure ${\cal
J}_2$ is never integrable. The almost Hermitian
structure $({\cal J}_2,h_t)$ is semi-K\"ahler. It is
indefinite almost K\"ahler  or nearly K\"ahler if and
only if the metric $g$ is Einstein, self-dual and
$t\tau=12$ or $t\tau=-6$, respectively.
\end{enumerate}
\end{thm}

\pf
To see when ${\cal J}_1$ is integrable, let us note
first that the vertical space at any point
${\sigma}\in{\cal Z}$ is spanned by the vectors of the
form $V=X\wedge Y-K_{\sigma}X\wedge K_{\sigma}Y,
X,Y\in T_{\pi(\sigma})M$. Moreover, if $V$ is of this
form, then $\sigma\times V= K_{\sigma}X\wedge
Y+X\wedge K_{\sigma}Y$. Therefore, by Corollary 1, the
Nijenhuis tensor of ${\cal J}_1$ vanishes if and only
if
$$
R(V)\sigma + \sigma\times R(\sigma\times V)\sigma=0
$$
for every $\sigma\in{\cal Z}$ and $V\in{\cal
V}_{\sigma}$. In view of (3.2), this is equivalent to
$$
g({\cal R}(V),W)=g({\cal R}(\sigma\times
V),\sigma\times W)
$$
for every $\sigma\in{\cal Z}$, $V,W\in{\cal
V}_{\sigma}$. Now, varying
$\sigma=y_1s_1+y_2s_2+y_3s_3$ on the fibre
$y_1^2-y_2^2-y_3^2=1$ of ${\cal Z}$ over a point $p\in
M$, we see that the latter condition is
satisfied if and only if $g({\cal
R}(s_1),s_1)=-g({\cal R}(s_2),s_2)=-g({\cal
R}(s_3),s_3)$ and
$g({\cal R}(s_i),s_j)=0$ for $i\neq j$. These
identities are equivalent to the self-duality of the
metric $g$.

  The second identity of Corollary 1 shows that the
almost complex structure ${\cal J}_2$ is never
integrable.

  Corollary 3 and (3.2) imply that $({\cal J}_k,h_t)$,
$k=1,2$, is semi-K\"ahler (i.e.
$\delta\Omega_{k,t}=0$)
if and only if $g({\cal W}^-(\sigma),\sigma\times
U)=0$ for every $\sigma\in{\cal Z}$ and $U\in{\cal
V}_{\sigma}$ which is equivalent to ${\cal W}^-=0$.

  It follows from Corollary 2 that the fundamental
$2$-form of the
almost Hermitian structure $({\cal J}_k,h_t)$ is
closed if and only if for any $\sigma\in{\cal Z}$ and
$V\in{\cal V}_{\sigma}$ we have
$$
(-1)^kt{\cal R}(V)-2V=0
$$
This is equivalent to $g$ being Einstein, self-dual
metric with $t\tau=12(-1)^k$. In this case
the structure $({\cal J}_1,h_t)$ is indefinite
K\"ahler since the almost complex structure
${\cal J}_1$ is integrable.




   If the structure $({\cal J}_2,h_t)$ is nearly
K\"ahler, then
$$
(\nabla_{X^h}\Omega_{2,t})(Y^h,V)_{\sigma}+(\nabla_{Y^h}\Omega_{2,t})(X^h,V)_{\sigma}=0
$$
for every $\sigma\in{\cal Z}$, $V\in {\cal
V}_{\sigma}$, $X,Y\in T_{\pi(\sigma)}M$. This identity
and Lemma 2 imply
$$ g({\cal R}(\sigma\times V),X\wedge
K_{\sigma}Y-K_{\sigma}X\wedge Y)=0 $$
Now, taking into account Lemma 2, (3.3) and the latter
equality, we obtain
$$
0=(\nabla_{V}\Omega_{2,t})(X^h,Y^h)_{\sigma}-(\nabla_{X^h}\Omega_{2,t})(Y^h,V)_{\sigma}=
$$
$$-g(V,X\wedge Y-K_{\sigma}Y\wedge
K_{\sigma}X)-{t\over 4}g({\cal R}(V),X\wedge
Y-K_{\sigma}Y\wedge K_{\sigma}X)$$
$$-{{3t}\over 4}g({\cal R}(\sigma\times V),X\wedge
K_{\sigma}Y+K_{\sigma}X\wedge Y)\eqno
(3.6)$$
As we have mentioned, the vectors of the form
$W=X\wedge Y-K_{\sigma}X\wedge K_{\sigma}Y$, $X,Y\in
T_{\pi({\sigma})}M$, span the vertical space ${\cal
V}_{\sigma}$ and $\sigma\times W=X\wedge K_{\sigma}Y+
K_{\sigma}X\wedge Y$. Therefore (3.6) is equivalent to

$$
g(V,W)+{t\over 4}g({\cal R}(V),W)+{{3t}\over 4}g({\cal
R}(\sigma\times V),\sigma\times W)=0
$$
for every $V,W\in{\cal V}_{\sigma}$. Varying $\sigma$
on the fibres of ${\cal Z}$ we see that $g$ is
Einstein and self-dual, and that
$t\tau=-6$.

    Conversely, it is not hard to show that under these
conditions the structure $({\cal J}_k,h_t)$ on the
twistor space ${\cal Z}$ is nearly-K\"ahler.
\qed

\vskip 20pt

\noindent
{\bf 4. Holomorphic functions}
\vskip 5pt
\noindent
On the classical twistor space over a Riemannian
4-manifold with either almost
  complex structure, there are no global non-constant
holomorphic functions, even
  when the base manifold is
non-compact
\cite{DM, DMG}.  However for the hyperbolic twistor
space there is
  considerable difference from the classical case as we
shall see.

First we remark that for local existence of
holomorphic functions, the
  situation is the same in both the classical and
hyperbolic cases.  A $C^\infty$ function on an almost
complex manifold is said to be {\it holomorphic} if
its differential is complex-linear with respect to the
almost complex structure.  On the twistor spaces, for
any $n=0,1,2,3$, let ${\cal F}_n({\cal J}_i)$
denote the (possibly empty) set of points $\sigma$
such that $n$ is the maximal number of local ${\cal
J}_i$-holomorphic
functions with ${\Bbb C}$-linearly independent
differentials at $\sigma$.  In \cite{DM} it is shown
that  for the
classical twistor space ${\cal Z}$ of an oriented
Riemannian 4-manifold $M$ we have
${\cal Z}={\cal F}_0({\cal J}_1)\cup{\cal F}_3({\cal
J}_1)
={\cal F}_0({\cal J}_2)\cup{\cal F}_1({\cal J}_2)$;
moreover
$${\cal F}_3({\cal J}_1)=\pi^{-1}(Int\{p\in M: {\cal
W}^{-}_p=0\}),$$
$${\cal F}_1({\cal J}_2)=\pi^{-1}(Int\{p\in M: {\cal
R}_p={\cal W}^{+}_p\})$$
The same arguments
give this result for the hyperbolic twistor space as
well.

    Let $M$ be a pseudo-Riemannian four-manifold with
metric $g$ of signature $(2,2)$.
We shall say that $M$ is {\it neutral almost
hyperhermitian} if it admits a globally defined almost
quaternionic structure of the second kind
$(J_1,J_2,J_3)$ compatible with the metric $g$. If the
structure tensors $J_1,J_2,J_3$ are integrable (i.e.
their Nijenhuis tensors vanish) the manifold is called
a {\it neutral hyperhermitian} surface; if, moreover,
$J_1,J_2,J_3$ are parallel with respect to the
Levi-Civita connection of $g$, $M$ is called a {\it
neutral hyperk\"ahler} surface.

    Suppose $(M,g,J_1,J_2,J_3)$ is a neutral almost
hyperhermitian four-manifold. Then all almost complex
structures $J_y=y_1J_1+y_2J_2+y_3J_3$, $y\in
H=\{(y_1,y_2,y_3) \\
\in {\Bbb R}^3:y_1^2-y_2^2-y_3^2=1\}$, are compatible
with the metric $g$ and determine the same orientation
on $M$. We shall always consider $M$ with this
orientation. As in the hyperhermitian case, if $M$ is
a neutral hyperhermitian surface, the metric $g$ is
self-dual and every almost complex structure $J_y$,
$y\in H$, is integrable. If, moreover, $M$ is neutral
hyperk\"ahler, then it is indefinite K\"ahler and
Ricci flat. It has been observed by Kamada \cite{Ka}
(see also \cite{P}) that any compact neutral
hyperk\"ahler surface is biholomorphic to a compex
torus or a primary Kodaira surface. He has also
obtained a description of all neutral hyperk\"ahler
structures on the latter surfaces.

     Given a neutral almost hyperhermitian
four-manifold $(M,g,J_1,J_2,J_3)$, denote by
$\pi:{\cal Z}\to M$
the hyperbolic twistor space of $(M,g)$. The
$2$-vectors corresponding to $J_1,J_2,J_3$ via (3.1)
form a global frame of $\bigwedge^-$ and we have a
natural projection $p:{\cal Z}\to H$ defined by
$p(\sigma)=(y_1,y_2,y_3)$ where
$K_{\sigma}=y_1J_1(x)+y_2J_2(x)+y_3J_3(x)$,
$x=\pi(\sigma)$. Thus ${\cal Z}$ is diffeomorphic to
$M\times H$ by the map $\sigma\to
(\pi(\sigma),p(\sigma))$. Further, we shall consider
the hyperboloid $H$ with the complex structure $S$
determined by the restriction to $H$ of the metric
$-dy_1^2+dy_2^2+dy_3^2$ of ${\Bbb R}^3$, i.e.
$SV=y\times V$ for $V\in T_yH$ where
$\times$ is the vector product on ${\Bbb R}^3$ defined
by means of the paraquaternionic algebra.
It is obvious that $p$ maps any fibre of ${\cal Z}$
biholomorphically on $H$ with respect to ${\cal J}_1$
and $S$.

    The hyperboloid $H$ has two connected component
$H^{\pm}=\{(y_1,y_2,y_3)
\in {\Bbb R}^3:y_1^2-y_2^2-y_3^2=1, \pm y_1>0\}$ and
the antipodal map $y\to -y$ sends $H^+$
anti-holomorphically onto $H^-$. Note also that the
stereographic projection
$\displaystyle{(y_1,y_2,y_3)\to {{y_2+iy_3}\over {1\pm
y_1}}}$ of $H^{\pm}$ from the point
$(\mp 1,0,0)$ is $\pm$-biholomorphic onto the unit
disk $\Delta$ in the complex plane ${\Bbb C}$.

\begin{thm}
    Let $M$ be a neutral almost hyperhermitian
four-manifold and ${\cal Z}$ its hyperbolic twistor
space. Then the natural projection $p:{\cal Z}\to H$
is ${\cal J}_1$-holomorphic (resp. ${\cal
J}_2$-anti-holomorphic) if and only if $M$ is neutral
hyperhermitian (resp. neutral hyperk\"ahler).
\end{thm}

\pf
  Let $(g,J_1,J_2,J_3)$ be the neutral almost
hyperhermitian structure on $M$.

  As we have already mentioned, the restriction of $p$
to any fibre of ${\cal Z}$ is ${\cal
J}_1$-holomorphic. Therefore $p$ is ${\cal
J}_1$-holomorphic on ${\cal Z}$ if and only if
$p_{\ast}((K_{\sigma}X)^h_{\sigma})=Sp_{\ast}(X^h_{\sigma})$
for any $\sigma\in {\cal Z}$ and
$X\in T_{\pi(\sigma)}M$. Given $\sigma\in {\cal Z}$,
the endomorphism $K_{\sigma}$ has the form
$K_{\sigma}=y_1J_1(x)+y_2J_2(x)+y_3J_3(x)$,
$x=\pi(\sigma)$, where $y=(y_1,y_2,y_3)\in H$ and we
set $J_y=y_1J_1+y_2J_2+y_3J_3$. Then, by (2.1), we
have
$$p_{\ast}(X^h_{\sigma})=-\sum_{b=1}^3
\epsilon_b(\langle D_XJ_y,J_b\rangle\circ\pi)
{\partial\over\partial y_b}$$
and it follows that $p$ is ${\cal J}_1$ holomorphic if
and only if
$$D_{J_yX}J_y=J_yD_{X}J_y$$
for any $y\in H$ and $X\in TM$. The latter condition
is equivalent to the almost complex structures
$J_y,y\in H,$ being integrable.

     Similarly, the projection $p$ is ${\cal
J}_2$-anti-holomorphic if and only if
$$D_{J_yX}J_y=-J_yD_{X}J_y$$
for any $y\in H$ and $X\in TM$. The latter condition
is equivalent to the almost complex structures $J_y$
being quasi K\"ahler. In dimension four this is
equivalent to $J_y$ being almost K\"ahler (i.e. with
closed fundamental $2$-forms).

    Now the theorem follows from the following:
\begin{lem}

\noindent
\begin{enumerate}
\item[$(i)$] The almost complex structures $J_y$,
$y\in H$, are integrable if and only if $M$ is neutral
hyperhermitian;
\item[$(ii)$] The almost complex structures $J_y$,
$y\in H$, are almost K\"ahler if and only if $M$ is
neutral hyperk\"ahler.
\end{enumerate}
\end{lem}

{\bf Proof of the lemma:} The almost complex structure
$J_y$ is integrable if and only if
$\displaystyle{D_{J_yX}J_y=J_yD_{X}J_y}$ for any $X\in
TM$. This identity is fulfilled for every $y\in H$
if and only if
$$
D_{J_kX}J_l+D_{J_lX}J_k=J_kD_{X}J_l+J_lD_{X}J_k, \>
1\leq k,l\leq 3.
$$
Using (2.6) and the paraquaternionic identities for
$J_1,J_2,J_3$, we see that, in the notation of (2.6),
the latter equalities are satisfied if and only if
$\alpha(X)=-\beta(J_1X)=-\gamma(J_2X)$ which is
equivalent to the integrability of the structure
tensors $J_1,J_2,J_3$.

    Proceeding in the same way, we see that
$\displaystyle{D_{J_yX}J_y=-J_yD_{X}J_y}$ for any
$y\in H$ and $X\in TM$ if and only if
$\alpha=\beta=\gamma=0$, i.e. $DJ_1=DJ_2=DJ_3=0$.
\qed

\begin{cor}
   Let $M$ be a compact neutral hyperhermitian manifold
with hyperbolic twistor space ${\cal Z}$ and let
$p:{\cal Z}\to H$ be the natural projection. Then any
${\cal J}_1$-holomorphic function $f$ on ${\cal Z}$
has the form $f=g\circ p$ where $g$ is a holomorphic
function on $H$. If $M$ is neutral hyperk\"ahler, any
${\cal J}_2$-holomorphic function $f$ on ${\cal Z}$
has the form $f=g\circ p$ where $g$ is an
anti-holomorphic function on $H$.
\end{cor}

\pf Any global section $s:N\to {\cal Z}$ of the
hyperbolic twistor bundle of an oriented
pseudo-Riemannian four-manifold $N$ with a neutral
metric determines a compatible almost complex
structure $K_s$ on $N$ via (3.1) and vice versa. Since
$s_{\ast}(X)=X^h\circ s+D_Xs$ for any $X\in TN$, it
follows from (3.3) that the map $s:(N,K_s)\to ({\cal
Z},{\cal J}_1)$ is holomorphic if and only if $K_s$ is
integrable; $s:(N,K_s)\to ({\cal Z},{\cal J}_2)$ is
holomorphic if and only if $K_s$ is almost K\"ahler.

    Now let $f$ be a ${\cal J}_1$-holomorphic function
on the twistor space ${\cal Z}$ of $M$. For any $y\in
H$, denote by $s_y$ the section of ${\cal Z}$
determined by the almost complex structure
$J_y=y_1J_1+y_2J_2+y_3J_3$. By Lemma 3(i), the
structure $J_y$ is integrable, hence the map
$s_y:(M,J_y)\to ({\cal Z},{\cal J}_1)$ is holomorphic.
Therefore $f\circ s_y$ is a holomorphic function on
the compact manifold $M$. So $f\circ s_y$ is a
constant and defining a function $g$ on $H$ by
$g(y)=f\circ s_y$ we have $f=g\circ p$. Since the
restriction of $p$ on a fibre of ${\cal Z}$ is a
biholomorphism onto $H$, the function $g$ is
holomorphic. Conversely, if $g$ is a holomorphic
function on $H$, then $f=g\circ p$ is a ${\cal
J}_1$-holomorphic function on ${\cal Z}$ by Theorem 3.

    Similar arguments prove the statement for the
${\cal J}_2$-holomorphic functions on ${\cal Z}$.
\qed

    In \cite{P} J.Petean has classified the compact
complex surfaces that admit
indefinite K\"ahler-Einstein metrics. In particular,
he has explicitly constructed
Ricci flat (non-flat) examples of such metrics on
complex tori,
hyperelliptic surfaces and primary Kodaira surfaces.

     Next we shall examine the J.Petean metrics on
${\Bbb R}^4$, a non-compact
manifold. All of them have the form:
$$g=f(dx_1\otimes dx_1+dx_2\otimes dx_2)+ dx_1\otimes
dx_3+ dx_3\otimes dx_1 +
dx_2\otimes dx_4 + dx_4\otimes dx_2$$
where $(x_1,x_2,x_3,x_4)$ are the standard coordinates
on ${\Bbb R}^4$ and
$f$ is a smooth positive function depending on $x_1$
and $x_2$ only.
Consider the frame given by
$$
{\bf e}_1=\frac{1}{\sqrt f}\frac{\partial}{\partial
x_1},\quad
{\bf e}_2=\frac{1}{\sqrt f}\frac{\partial}{\partial
x_2},$$
$${\bf e}_3=-\frac{1}{\sqrt f}\frac{\partial}{\partial
x_1}+
\sqrt f\frac{\partial}{\partial x_3},\quad
{\bf e}_4=-\frac{1}{\sqrt f}\frac{\partial}{\partial
x_2}+
\sqrt f\frac{\partial}{\partial x_4}
$$
Then $||{\bf e}_1||=||{\bf e}_2||=-||{\bf
e}_3||=-||{\bf e}_4||=1$ and let $s_1,s_2,s_3$ be the
sections of $\Lambda^{-}{\Bbb R}^4$ defined as in the
last section.
A direct computation shows that these sections are
parallel, therefore they
define a neutral hyperk\"ahler structure on ${\Bbb
R}^4$. So, the
metric $g$ is Ricci flat and self-dual \cite{P}, hence
the almost complex structure
${\cal J}_1$ on the hyperbolic twistor space ${\cal
Z}$ of $({\Bbb R}^4,g)$ is integrable. Moreover,
we have the following result which shows that there
can be an abundance of holomorphic functions
on a hyperbolic twistor space.
\begin{thm}
The hyperbolic twistor space $({\cal Z},{\cal J}_1)$
of $({\Bbb R}^4,g)$ is
biholomorphic to ${\Bbb C}^2\times H$.
\end{thm}

\pf
   Since the sections $s_1,s_2,s_3$ of ${\cal Z}$ are
globally defined, ${\cal Z}$ is diffeomorphic to
${\Bbb R}^4\times H$. Denote by ${\cal Z}^{\pm}$ the
connected component of ${\cal Z}$ determined by the
hyperbolic plane $H^{\pm}=\{(y_1,y_2,y_3)\in {\Bbb
R}^3: y_1^2-y_2^2-y_3^2=1,\pm  y_1>0\}$. We shall
identify $H^{\pm}$ with the unit disk $\Delta$ in the
complex plane ${\Bbb C}$ by means of the
``stereographic" projection
$\displaystyle{(y_1,y_2,y_3)\to {{y_2\pm iy_3}\over
{1\pm y_1}}}$ and ${\cal Z}^{\pm}$ with ${\Bbb
R}^4\times\Delta$.   Our proof will be to show that
there are global ${\cal J}_1$-holomorphic coordinates
on ${\cal Z^{\pm}}\cong{\Bbb R}^4\times\Delta$. We
shall consider only the component
${\cal Z}^+$ of ${\cal Z}$ since the same reasoning
works for ${\cal Z}^-$.

    Since the sections $s_1,s_2,s_3$ are parallel, it
follows from (2.1) that the complex structure ${\cal
J}_1$ on ${\Bbb R}^4\times\Delta$ is given by
$$
{\cal J}_1{\bf e}_1=y_1{\bf e}_2+y_2{\bf e}_3+y_3{\bf
e}_4,\quad
{\cal J}_1{\bf e}_2=-y_1{\bf e}_1+y_3{\bf e}_3-y_2{\bf
e}_4,
$$
$$
{\cal J}_1{\bf e}_3=y_2{\bf e}_1+y_3{\bf e}_2+y_1{\bf
e}_4,\quad
{\cal J}_1{\bf e}_4=y_3{\bf e}_1-y_2{\bf e}_2-y_1{\bf
e}_3,
$$
$$
{\cal J}_1\frac{\partial}{\partial
x}=\frac{\partial}{\partial y},\quad
{\cal J}_1\frac{\partial}{\partial
y}=-\frac{\partial}{\partial x}
$$
where $x,y$ are the standard coordinates on $\Delta$
and $y_1,y_2,y_3$ are defined by
$z=x+iy$ as follows
$$
y_1=\frac{1+|z|^2}{1-|z|^2},\quad y_2=\frac{z+\bar
z}{1-|z|^2},\quad
y_3=\frac{z-\bar z}{i(1-|z|^2)}.
$$
Note that
$${\bf e}_1+i{\cal J}_1{\bf e}_1, \quad {\bf
e}_3+i{\cal J}_1{\bf e}_3,\quad
\displaystyle{\frac{\partial}{\partial x}+{\cal
J}_1\frac{\partial}{\partial x}}$$
is a global frame of the bundle $T^{0,1}({\Bbb
R}^4\times\Delta)$ of $(0,1)$-vectors with respect to
${\cal J}_1$. Thus a smooth complex-valued function
$G$ on ${\Bbb R}^4\times\Delta$ is ${\cal
J}_1$-holomorphic if and only if it satisfies the
Cauchy-Riemann equations $$({\bf e}_1+i{\cal J}_1{\bf
e}_1)G=({\bf e}_3+i{\cal J}_1{\bf e}_3)G=
\displaystyle{(\frac{\partial}{\partial x}+{\cal
J}_1\frac{\partial}{\partial x})}G
=0.$$

    Clearly the projection
$G_1:{\Bbb R}^4\times\Delta\to\Delta$ is a ${\cal
J}_1$-holomorphic function, therefore any ${\cal
J}_1$-holomorphic function
$G$ on ${\Bbb R}^4\times\Delta$ can be expended in a
series of $z=x+iy$ with coefficients being smooth
functions on ${\Bbb R}^4$. This remark leads us to
seek a ${\cal J}_1$-holomorphic function which is
linear in $z$ and
it is easy to check that the function
$G_2=(x_1+ix_2)+iz(x_1-ix_2)$ is ${\cal
J}_1$-holomorphic. Next we shall show that there
exists a third ${\cal J}_1$-holomorphic function $G_3$
on ${\Bbb R}^4\times\Delta$ such that
$G_1, G_2, G_3$ form global ${\cal J}_1$-holomorphic
coordinates on ${\Bbb R}^4\times\Delta$.
To do this, we take $G_1$ and $G_2$ as new
coordinates, i.e. we introduce new smooth
coordinates on ${\Bbb R}^4\times\Delta$ by setting
$$p=x_1(1-y)+x_2x,\; q=x_1x+x_2(1+y),\; r=x_3,\;
s=x_4,\;u=x,\;v=y \eqno(4.1)$$
Set
$$F(p,q,u,v)=f(\frac{(1+v)p-uq}{1-u^2-v^2},
\frac{-up+(1-v)q}{1-u^2-v^2})=f(x_1,x_2)\eqno(4.2)$$
It is straightforward to compute the action of ${\cal
J}_1$ in the new coordinates (4.1) and to see that
$$
  \frac{\partial}{\partial p}+i{\cal
J}_1\frac{\partial}{\partial p}=
\frac{\partial}{\partial p}+i\frac{\partial}{\partial
q}+i\frac{2uF}{(1-u^2-v^2)^2}\frac{\partial}{\partial
r}+i\frac{2(u^2+v^2+v)}{(1-u^2-v^2)^2}\frac{\partial}{\partial
s}
$$
$$
\frac{\partial}{\partial r}+i{\cal
J}_1\frac{\partial}{\partial r}=
(1+i\frac{2u}{1-u^2-v^2})\frac{\partial}{\partial
r}+i\frac{1+u^2+v^2+2v}{1-u^2-v^2}\frac{\partial}{\partial
s}
$$
$$
\frac{\partial}{\partial u}+i{\cal
J}_1\frac{\partial}{\partial u}=
\frac{\partial}{\partial u}+i\frac{\partial}{\partial
v}+
i\frac{2[(2u^2+u^2v+v^3-v)p-u(1+u^2+v^2-2v)q]F}{(1-u^2-v^2)^3}\frac{\partial}{\partial
r}
$$
$$
+i\frac{2[u(1+u^2+v^2+2v)p-(2u^2-u^2v-v^3+v)q]F)}{(1-u^2-v^2)^3}\frac{\partial}{\partial
s}
$$
The vector fields
$$
\frac{\partial}{\partial p}+i{\cal
J}_1\frac{\partial}{\partial p},\quad
\frac{\partial}{\partial r}+i{\cal
J}_1\frac{\partial}{\partial r},\quad
\frac{\partial}{\partial u}+i{\cal
J}_1\frac{\partial}{\partial u}
$$
form a global frame for the bundle $T^{0,1}({\Bbb
R}^4\times\Delta)$ of $(0,1)$-vectors. Thus a smooth
complex-valued function $G(p,q,r,s,u,v)$
on ${\Bbb R}^4\times\Delta$ is ${\cal
J}_1$-holomorphic if and only if the Cauchy-Riemann
equations
$$
(\frac{\partial}{\partial p}+i{\cal
J}_1\frac{\partial}{\partial
p})G=(\frac{\partial}{\partial r}+i{\cal
J}_1\frac{\partial}{\partial
r})G=(\frac{\partial}{\partial u}+i{\cal
J}_1\frac{\partial}{\partial u})G=0
$$
are satisfied. Now set $z=u+iv$ and $w=p+iq$.
Then a direct computation
shows that $G$ is ${\cal J}_1$-holomorphic if and only
if
$$\frac{\partial G}{\partial
s}=-i\frac{z-i}{z+i}\frac{\partial G}{\partial
r}\eqno(4.3)$$
$$\frac{\partial G}{\partial\overline w}=
\frac{zF}{(1-|z|^2)(z+i)}\frac{\partial G}{\partial
r}\eqno(4.4)$$
$$\frac{\partial G}{\partial\overline z}=
\frac{(iw+z\overline
w)zF}{(1-|z|^2)(z+i)}\frac{\partial G}{\partial
r}\eqno(4.5)$$
As we have mentioned the functions
$$G_1(p,q,r,s,u,v)=u+iv=z,\quad
G_2(p,q,r,s,u,v)=p+iq=w$$
are ${\cal J}_1$-holomorphic.
We shall seek a third holomorphic function $G_3$ in
the form
$$G_3(p,q,r,s,u,v)=r-i\frac{z-i}{z+i}s+H(p,q,u,v)\eqno(4.6)$$
where $H$ is a smooth function on ${\Bbb
R}^2\times\Delta$. The function $G_3$ satisfies (4.3),
(4.4) and (4.5) provided
$$\frac{\partial H}{\partial\overline w}=
\frac{zF}{(1-|z|^2)(z+i)},\quad
\frac{\partial H}{\partial\overline z}=
\frac{(iw+z\overline
w)zF}{(1-|z|^2)^2(z+i)}\eqno(4.7)$$
This system (which is, in fact, a
$\bar\partial$-equation on ${\Bbb
C}\times\Delta$)
has a global solution if and only if
$$
\frac{\partial}{\partial\overline
z}\Big(\frac{zF}{1-|z|^2)(z+i)}\Big )=
\frac{\partial}{\partial\overline
w}\Big(\frac{(iw+z\overline
w)zF}{(1-|z|^2)(z+i)}\Big )
$$
which (in view of (4.2)) is equivalent to the
following identity:
$$\frac{\partial f}{\partial x_1}.\frac{\partial x_1}
{\partial\overline z}+
\frac{\partial f}{\partial x_2}.\frac{\partial
x_2}{\partial\overline z}
=
\frac{iw+z\overline w}{(1-|z|^2)}
\Big(\frac{\partial f}{\partial x_1}.\frac{\partial
x_1}{\partial\overline w}+
\frac{\partial f}{\partial x_2}.\frac{\partial
x_2}{\partial\overline w}\Big )\eqno(4.8)$$
On the other hand
$$x_1=\frac{(2i+z-\overline z)(w+\overline
w)-(z+\overline z)(w-\overline w)}
{4i(1-|z|^2)}$$
and
$$x_2=- \frac{(z+\overline z)(w+\overline
w)+(2i-z+\overline z)(w-\overline w)}
{4(1-|z|^2)}.$$
Hence
$$\frac{\partial x_1}{\partial\overline z}=
\frac{(z+i)(w-iz\overline w)}{2(1-|z|^2)^2},\quad
\frac{\partial x_2}{\partial\overline z}=
-i\frac{(z-i)(w-iz\overline w)}{2(1-|z|^2)^2}$$
$$\frac{\partial x_1}{\partial\overline w}=
-i\frac{z+i}{2(1-|z|^2)},\quad
\frac{\partial x_2}{\partial\overline w}=
-\frac{z-i}{2(1-|z|^2)}.$$
These identities show that
$$
\frac{\partial x_1}{\partial\overline z}=
\frac{iw+z\overline w}{1-|z|^2}
\frac{\partial x_1}{\partial\overline w},
\quad
\frac{\partial x_2}{\partial\overline z}=
\frac{iw+z\overline w}{1-|z|^2}
\frac{\partial x_2}{\partial\overline w}
$$
and therefore identity (4.8) is satisfied.

     Let $H(p,q,u,v)$ be a global solution of (4.7).
Consider the map
$(G_1,G_2,G_3):({\Bbb R}^4\times\Delta,J_1)\to {\Bbb
C}^2\times\Delta$.
It is holomorphic and bijective since for any
$(\alpha,\beta,\gamma)\in
{\Bbb C}^2\times\Delta$ the system $u+iv=\alpha,
p+iq=\beta,
r-\displaystyle{i\frac{z-i}{z+i}s}+H(p,q,u,v)=\gamma$
has a unique solution ($u, v$ and
$p, q$ are uniquely
determined by $\alpha$ and $\beta$, then $r, s$  are
uniquely determined by
$\gamma-H(p,q,u,v)$ since
$\displaystyle{\overline{i\frac{z-i}{z+i}}\neq
i\frac{z-i}{z+i}}$ for every
$z$ with $|z|\neq 1$).
\qed

\noindent
{\bf Remark}. The choice of the function $G_3$ is
almost canonical since by the
first equation of (4.7) it has the form
$$
G_3=H(u,v,p,q,r-i\frac{z-i}{z+i}s)
$$
and it can be shown that $H$ is a holomorphic function
of
$t=r-\displaystyle{i\frac{z-i}{z+i}s}$.
Then for any fixed $(u,v,p,q)$, $t\to G_3(u,v,p,q,t)$
is a biholomorphism
of ${\Bbb C}$ , hence a linear function of $t$.

\vskip 20pt
\noindent
Department of Mathematics\hskip1cm Institute of
Mathematics and Informatics

\noindent
Michigan State University\hskip1.3cm Bulgarian Academy
of Sciences

\noindent
East Lansing, MI 48824\hskip1.7cm Sofia, Bulgaria

\noindent
blair@math.msu.edu\hskip2.2cm jtd@math.bas.bg;
muskarov@math.bas.bg



\begin{thebibliography}{99}


\bibitem{AGI}
Alexandrov, B., Grantcharov, G. and Ivanov, S.:
Curvature properties of twistor
spaces of quaternionic K\"ahler manifolds, {\it J. Geom.}
{\bf 62} (1998), 1-12.

\bibitem{AHS}
Atiyah, M. F., Hitchin, N. J. and Singer I.M.:
Self-duality  in
four-dimensional Riemannian  geometry,
{\it Proc. Roy. Soc. London}  Ser. A {\bf 362} (1978), 425-461.

\bibitem{bes}
Besse A. L.: {\it Einstein Manifolds}, Springer, Berlin,
1987.

\bibitem{DB}
Blair D. E.: A hyperbolic twistor space, {\it Balkan
J. Geom. and Appl.} {\bf 5} (2000),
9-16.

\bibitem{Bl}
Bla\u zi\'c N., Paraquaternionic projective space
and
pseudo-Riemannian geometry, {\it  Pub. Inst. Math.} {\bf 60}
(1996), 101-107.

\bibitem{DM}
Davidov, J. and Mu\u skarov, O.: Existence of
holomorphic functions on twistor spaces,
{\it  Bull. Soc. Math. Belgique} {\bf 40} Ser. B (1989),
131-151.

\bibitem{DMG}
Davidov, J. and Mu\u skarov, O. and Grantcharov, G.:
Almost complex
structures on twistor spaces, {\it Almost Complex
Structures}, World
Scientific, Singapore, 1994, 113-149.

\bibitem{ES}
Eells, J. and Salamon, S.: Twistorial construction
of harmonic maps
of surfaces into four-manifolds, {\it Ann. Scuola Norm.
Sup.  Pisa}  Cl.\,Sci.
{\bf 12} (1985), 589-640.

\bibitem{GMV}
Garcia-Rio, E., Matsushita,Y. and V\'azquez-Lorenzo, R.:
  Paraquaternionic K\"ahler manifolds, {\it Rocky Mountain J. Math.} {\bf
31} (2001), 237-260.

\bibitem{G}
Gauduchon, P.: Structures de Weyl et th\'eor\`ems
d'annualation sur
une vari\'et\'e conforme autoduale,  {\it Ann.Scuola
Norm.Sup.}, ser.IV {\bf 18} (1991),
563-629.

\bibitem{Ia}
Ianus, S.: Sulle strutture canoniche dello spazio
fibrato tangente di
una variet\`a riemanniana, {\it Rend. Mat.} {\bf 6} (1973),
1-22.

\bibitem{IU}
Ianus, S. and Udriste, C.: Asupra spatiului fibrat
tangent al unei
varietati diferentiabile, {\it St. Cerc. Mat.} {\bf 22}
(1970), 599-611.

\bibitem{Is}
Ishihara, S.: Quaterion K\"ahlerian manifolds, {\it J.
Differential Geometry} {\bf 9} (1974), 483-500.

\bibitem{Ka}
Kamada, H.: Neutral hyperk\"ahler structures on
primary Kodaira
surfaces, {\it Tsukuba J. Math.} {\bf 23} (1999), 321-332.

\bibitem{Lib}
Libermann, P.: Sur les structures presque
quaternioniennes de
deuxi\`eme esp\`ece, {\it C. R. Acad. Sci Paris} {\bf 234}
(1952), 1030-1032.

\bibitem{M}
Mu\u skarov, O.: Almost Hermitian structures on
twistor spaces and
their type, {\it Atti Sem.Mat.Fis.Univ.Modena} {\bf 37}
(1989), 285-297.

\bibitem{P}
Petean, J. Indefinite K\"ahler-Einstein metrics
on compact complex surfaces,
{\it Commun. Math. Phys.} {\bf 189} (1997), 227-235.

\bibitem{Sal}
Salamon, S.: Quaternionic K\"ahler manifolds, {\it Invent.
Math.} {\bf 67} (1982),
143-171.

\end{thebibliography}
\end{document}